\documentclass[10pt]{article}

\newcounter{mascotsection}
\newcounter{mascotsubsection}[mascotsection]

\newcommand{\mascotsection}[1]{%
	\stepcounter{mascotsection}
	\section{\themascotsection . #1}
}
\newcommand{\mascotsubsection}[1]{%
	\stepcounter{mascotsubsection}
	\subsection{\themascotsection .\themascotsubsection . #1}
}

\usepackage{epsfig}

\usepackage{amsfonts,amsmath,amssymb,amsthm}
\usepackage{graphicx}

\newcommand{\RR}{\mathbb{R}}

\begin{document}
\vskip 1.8cm 
\centerline{\LARGE\bf Hybrid spherical approximation}
\vskip 0.6cm
\centerline{\bf Alessandra De Rossi}
\centerline{Department of Mathematics \lq\lq G. Peano\rq\rq, University of Torino,}
\centerline{Via Carlo Alberto 10,~10123 Torino,~Italy}
\centerline{\em alessandra.derossi@unito.it}
\vskip .6cm
\begin{abstract}
\noindent In this paper a local approximation method on the sphere is presented. As interpolation scheme we consider a partition of unity method, such as the modified spherical Shepard's method, which uses zonal basis functions (ZBFs) plus spherical harmonics as local approximants. Moreover, a spherical zone algorithm is efficiently implemented, which works well also when the amount of data is very large, since it is based on an optimized searching procedure. Numerical results show good accuracy of the method, also on real geomagnetic data. 
\end{abstract}

\medskip

\medskip
\mascotsection{Introduction}

Over the last decades approximation of functions on the sphere has attracted the interest of many researchers. In particular, the use of zonal basis functions (ZBFs) and spherical harmonics appears in a wide field of applications in numerical mathematics and computer science. Applications can be found in approximation of  scattered data, for example in geophysical and meteorological problems. These functions are of special interest, since they show several features which make them well suited for a wide range of problems and, at the same time, computationally attractive (see, e.g., \cite{Fasshauer-Schumaker98, Wendland05} and references therein).

In this paper, following the idea in \cite{Sloan-Sommariva08}, there analyzed in a global setting, we analyze a local interpolation scheme on the sphere, combining ZBFs with spherical harmonics of low degree. The aim of our paper is to verify if the addition of a polynomial part in a local approach, which is based on a classical partition of unity method as the well-known modified spherical Shepard's formula, allows us to improve accuracy of the considered interpolation technique. 
The basis of the spherical algorithm employed in the numerical experiments is the one presented and tested in \cite{Cavoretto-DeRossi10} (see also \cite{Cavoretto-DeRossi12}). It has been modified and efficiently updated for our purposes.

The paper is organized as follows. In Section 2 we consider some basic mathematical tools, referring to spherical harmonics and ZBFs. Section 3 is devoted to present the local Shepard's method which uses ZBFs phus spherical harmonics as local approximants. Section 4 refers to the spherical interpolation algorithm, while in Section 5 numerical results are presented.

\mascotsection{Functions on the sphere}

\mascotsubsection{Spherical harmonics}

We start this section by recalling the analogue of classical polynomials on the sphere, called \textsl{spherical harmonics} \cite{Fasshauer-Schumaker98}. Thus, denoting by ${\cal P}_d={\cal P}_d(\RR^3)$ the space of trivariate polynomials of degree at most $d$ and ${\cal H}_d={\cal H}_d(\mathbb{S}^{2})$ its restriction on the unit sphere, i.e. ${\cal H}_d:={\cal P}_d|_{\mathbb{S}^{2}}$, we know that a trivariate polynomial $p$ is called \textsl{homogeneous of degree $d$} provided $p(t\textbf{x})=t^d p(\textbf{x})$ for any $t > 0$ and any $x\in\RR^3$. It is called \textsl{harmonic} if $\Delta p(\textbf{x})\equiv0$, where $\Delta$ is the Laplace operator. Then, we can define the linear space of spherical harmonics of exact degree $d$ as follows: 
\begin{equation}
H_d:=\left\{p|_{\mathbb{S}^{2}}: p \in {\cal P}_d, \mbox{ $p$ homogeneous of degree $d$ and harmonic}\right\}.	\nonumber
\end{equation} 

Given the Laplace-Beltrami operator $\Delta_{2}$ on the unit sphere, the eigenvalues of the eigenvalue problem $(\Delta_{2}+\lambda)f = 0$ are $\lambda_d = d(d+1)$, $d\geq 0$, and the space $H_d$ is precisely the eigenspace of $\Delta_{2}$ correponding to $\lambda_d$. The dimension of $H_d$ is given by the multiplicity of $\lambda_d$, i.e. $N_{d,3}={\rm dim}(H_d)= 2d+1$ (see \cite{Hubbert-Morton04}).

It is known that $H_d$ is the orthogonal complement of ${\cal H}_{d-1}$ in the space ${\cal H}_d$ with respect to the $L_2$-inner product on $\mathbb{S}^{2}$
\begin{equation}
 (f,g)_{L_2(\mathbb{S}^{2})}=\int_{\mathbb{S}^{2}} f(\textbf{x})g(\textbf{x}) d\mu(\textbf{x}), \nonumber	
\end{equation}
where $d\mu(\textbf{x})$ is the standard measure on the sphere. Using this fact repeatedly, we have that
\begin{equation}
	{\cal H}_d = \oplus_{j=0}^{d} H_j. \nonumber
\end{equation}

Since the spherical harmonics form an orthonormal basis for $L_2(\mathbb{S}^{2})$, every function $f\in L_2(\mathbb{S}^{2})$ has a Fourier expansion. Thus, given an orthonormal basis ${\cal B}_d=\left\{Y_{d,k}: k=1, 2, \ldots, N_{d,3}\right\}$ for ${\cal H}_d$, the orthonormal system $\left\{{\cal B}_d\right\}_{d=0}^{\infty}$ is complete in $L_2(\mathbb{S}^{2})$, and every $f\in L_2(\mathbb{S}^{2})$ has a spherical Fourier representation of the form 
\begin{equation}
	\label{Fourier}
	f=\sum_{d=0}^{\infty}\sum_{k=1}^{N_{d,3}} \hat{f}_{d,k}Y_{d,k}, \nonumber
\end{equation}
where	$\hat{f}_{d,k}=\left(f, Y_{d,k}\right)_{L_2(\mathbb{S}^{2})}$ are the spherical Fourier coefficients of $f$. See e.g. \cite{Wendland05} for further details. 

\mascotsubsection{Zonal basis functions}

Let ${\cal \chi}_n=\{ (\textbf{x}_i,f_i), i=1, 2, \ldots, n\} \subset \mathbb{S}^{2} \times \RR$ be the set of pairs such that $\textbf{x}_i$ is a node and $f_i$ is the corresponding data value of an unknown function $f:\mathbb{S}^{2}\rightarrow \RR$. So the interpolation problem consists in finding a function $s:\mathbb{S}^{2}\rightarrow \RR$, which satisfies the interpolation conditions 
\begin{eqnarray}
 \label{condinterp} 
  s({\textbf{x}}_i)=f_i, \hspace{.5cm} i=1,2,\ldots,n.
\end{eqnarray}
The interpolating function $s$ might also be expressed as a linear combination of a {\sl zonal basis function} $\psi:[0,\pi]\rightarrow \RR$, i.e. 
\begin{eqnarray}
\label{intfun}
s({\textbf{x}})=\sum_{j=1}^n a_j \psi (g(\textbf{x},\textbf{x}_j)), \hskip0.5cm \textbf{x}\in \mathbb{S}^{2},
\end{eqnarray}
where $g(\textbf{x},\textbf{x}_j)=\arccos(\textbf{x}^T \textbf{x}_j)$ denotes the geodesic distance, and $s$ satisfies the interpolation conditions (\ref{condinterp}). Thus, we have uniqueness of the interpolation process if and only if the interpolation matrix $A \in \RR^{n \times n}$, which is given by
\begin{eqnarray}
\label{mat}
A_{i,j}=\psi(g(\textbf{x}_i,\textbf{x}_j)), \hskip 0.5cm 1\le i,j\le n,
\end{eqnarray}
turns out to be non-singular. In fact, even though there is no complete characterization for those functions satisfying the non-singularity condition, we know that a sufficient condition is that the matrix $A$ is positive definite (see \cite{Baxter-Hubbert01}). Moreover, the continuous function $\psi:[0,\pi] \rightarrow \RR$ is called {\sl positive definite} of order $n$ on $\mathbb{S}^{2}$, if, for any set of nodes,
\begin{eqnarray}
\label{posdef}
	\sum_{i=1}^n \sum_{j=1}^n a_i a_j \psi(g(\textbf{x}_i,\textbf{x}_j)) \geq 0,
\end{eqnarray}
for any $\textbf{a}=[a_1,a_2,\ldots,a_n]^T \in \RR^n$. The function $\psi$ is called \textsl{strictly positive definite} of order $n$ if the quadratic form (\ref{posdef}) is zero only for $\textbf{a}\equiv \textbf{0}$. If $\psi$ is strictly positive definite for any order $n$, then it is called \textsl{strictly positive definite}. Therefore, if $\psi$ is strictly positive definite, the interpolant (\ref{intfun}) is unique, the matrix (\ref{mat}) being positive definite and so non-singular.

Now, if we add a spherical harmonic of degree $d$ to a linear combination of the form (2), the interpolant takes the form
\begin{eqnarray}
\label{harm_interp}
s({\textbf{x}})=\sum_{j=1}^n a_j \psi (g(\textbf{x},\textbf{x}_j)) + \sum_{k=1}^U b_k Y_k (\textbf{x}), \hskip0.5cm \textbf{x}\in \mathbb{S}^{2},
\end{eqnarray}
where $U={\rm dim} {\cal H}_d(\mathbb{S}^{2})$, and $\{Y_1,Y_2,\ldots,Y_U\}$ is a basis for ${\cal H}_d (\mathbb{S}^{2})$.

The solution of the interpolation problem in the form given in (\ref{harm_interp}) is obtained by requiring that $s$ satisfies the interpolation conditions (\ref{condinterp}), and the additional conditions (see \cite{Fasshauer-Schumaker98})
\begin{eqnarray}
\label{addcond}
\sum_{i=1}^{n} a_i Y_k(\textbf{x}_i)=0, \hspace{.5cm} \hbox{for } k=1,2,\ldots,U. 
\end{eqnarray}
This problem consists in solving a system of $n$ linear equations in $n+U$ unknowns. Thus, assuming that $n \geq U$, we have the linear system
\begin{eqnarray}
\label{system}
\left[
\begin{array}{cc}
A   & Y  \\
Y^T & O 
\end{array}
\right]
\left[
\begin{array}{c}
\textbf{a}  \\
\textbf{b} 
\end{array}
\right]
=
\left[
\begin{array}{c}
\textbf{f}  \\
\textbf{0} 
\end{array}
\right],
\end{eqnarray}
where $A=\{\psi(g(\textbf{x}_i,\textbf{x}_j))\}$ is an $n\times n$ matrix (as in (\ref{mat})), $Y = \{Y_k(\textbf{x}_i)\}$ is an $n\times U$ matrix, and $f$ denotes the column vector of the function values $f_i$.


\mascotsection{Local Shepard's method} \label{mssm_sphere}
In this section we consider a modified version of spherical Shepard's method, which uses ZBFs plus spherical harmonics as local approximants. This approach exploits accuracy of ZBFs, overcoming some drawbacks such as instability and inefficiency of the global ZBF method.

Now, the modified spherical Shepard's interpolant $F:\mathbb{S}^{2} \rightarrow \RR$ is given by
\begin{eqnarray}
\label{modsphshep}
	F(\textbf{x})=\sum_{j=1}^{n} Z_j(\textbf{x}) \bar{W}_j(\textbf{x}), 
\end{eqnarray}
where
\begin{equation}
Z_j(\textbf{x})\equiv s_{\mid_{{\cal D}_j}}(\textbf{x}), \hspace{0.5cm} \textbf{x} \in {\cal D}_j\subset {\cal S}_n, \nonumber	
\end{equation}
is a local interpolant to $f$ in a vicinity of $\textbf{x}_j$, constructed on the restricted subset ${\cal D}_j$ containing the $n_Z$ nodes closest to $\textbf{x}_j$ and satisfying the interpolation conditions 
\begin{equation} \label{iii}
 Z_j(\textbf{x}_i)=f_i, \hspace{0.5cm} i=1,2,\ldots,n_Z.
\end{equation}
The weight functions $\bar{W}_j (\textbf{x})$, $j=1,2,\ldots,n$, are 
\begin{equation}
	\bar{W}_j(\textbf{x}) = \frac{W_j(\textbf{x})}{\sum_{k=1}^{n} W_k(\textbf{x})}, \hspace{0.5cm} j=1,2,\ldots,n, \nonumber
\end{equation}
with 
\begin{equation}
 W_j(\textbf{x}) = \tau(\textbf{x},\textbf{\textbf{x}}_j) / g(\textbf{x}, \textbf{x}_j). \nonumber	
\end{equation}
The localizing function $\tau(\textbf{x},\textbf{x}_j)$ is 
\begin{equation}
\tau(\textbf{x},\textbf{x}_j)= \left\{
\begin{array}{ll}
1, & \mbox{if $\textbf{x} \in {\cal C}(\textbf{x}_j;s)$}, \nonumber \\
0, & \mbox{otherwise}, 
\end{array}
\right.
\end{equation}
where ${\cal C}(\textbf{x}_j;r)$ is a spherical cap of centre at $\textbf{x}_j$ and spherical radius $r$. Note that the weights $\bar{W}_j$ constitute a partition of unity.

As regard to the choice of nodal functions we have a wide class of ZBFs, which are usually considered in the scattered data interpolation on the sphere. The nodal functions have the form 
\begin{eqnarray}
\label{zonbasfun}
  Z_j(\textbf{x}) = \sum_{i=1}^{n_Z}  a_i \psi(\arccos(\textbf{x}^T\textbf{x}_i))+\sum_{k=1}^U  b_k Y_k(\textbf{x}), \hspace{0.5cm} j=1,2,\ldots,n,
\end{eqnarray}
where the zonal basis functions $\psi(\arccos(\textbf{x}^T\textbf{x}_i))$ depend on the $n_Z$ nodes of the considered neighbourhood of $\textbf{x}_j$, and the space ${\cal H}_{v}$ spanned by the spherical harmonics $Y_k (\textbf{x})$ of degree $v$ has a dimension $U \leq n_Z$. Thus, we require that $Z_j$ satisfies the interpolation conditions (\ref{iii}) and the additional conditions 
\begin{equation}
\sum_{i=1}^{n_Z} a_i Y_k(\textbf{x}_i)=0, \hspace{.5cm} \hbox{for } k=1,2,\ldots,U. \nonumber
\end{equation}

We remark that, considering a strictly positive definite function $\psi$, we can generate ZBFs as the specialization on the sphere of the more general radial basis functions (RBFs). In fact, given any Euclidean RBF,
we may associate with it a ZBF \cite {Cavoretto-DeRossi10}. 
An example of strictly positive definite ZBF on $\mathbb{S}^{2}$ is the \textsl{spherical inverse multiquadric} (IMQ) \cite{Fasshauer-Schumaker98} 
\begin{equation}
\left.
\begin{array}{rcl}
\psi(t) & = & \left(1+\gamma^2 -2\gamma c\right)^{-1/2},
\end{array}
\right.
\end{equation}
where $\gamma\in \left(0,1\right)$, $c=\cos t$ and $t$ measures geodesic distance on the sphere.

\mascotsection{Spherical interpolation algorithm}
In this section we refer to the spherical algorithm, which is based on the partition of the sphere in spherical zones, that are portions of the spherical surface included between two parallel planes. The basis of this spherical interpolation algorithm has been proposed and widely tested in \cite{Cavoretto-DeRossi10}. Here, it has been modified adding spherical harmonics of low degree to the local ZBF interpolants. We remark that some details about the algorithm have been omitted, the interested readers can refer to \cite{Cavoretto-DeRossi10}.

For simplicity, we subdivide the description of the spherical algorithm in three parts, namely distribution, localization and evaluation phases.

\mascotsubsection{Distribution stage}
Let us consider the sets: ${\cal S}_n=\{(x_i,y_i,z_i), i=1,2,\ldots,n\}$, set of nodes, ${\cal F}_n=\{f_i, i=1,2,\ldots,n\}$ the set of the corresponding data values, and the set ${\cal E}_s=\{(x_i,y_i,z_i), i=1,2,\ldots,s\}$ of evaluation points. Then, we denote by $n_Z$ and $n_W$ the localization parameters.

Initially, the elements of the sets ${\cal S}_n$ and ${\cal E}_s$ are ordered with respect to the $z$-axis direction, applying a \textsl{quicksort procedure}. Then, for each node $(x_i,y_i,z_i)$, $i=1,2,\ldots,n$, a local circular neighbourhood (a spherical cap) is constructed, whose spherical radius depends on the node number $n$, the considered value $n_Z$, and the positive integer $k_1$ (which determines the radius of the spherical cap). More precisely, we define 
\begin{eqnarray}
\label{delta_Z}
 \delta_Z = \arccos \left(1-2\sqrt{k_1} \frac{n_Z}{n}\right), \hspace{0.5cm} k_1 = 1, 2, \ldots
 \end{eqnarray} 
After the number of spherical zones to be considered is found taking 
$q = \left\lceil \pi/\delta_Z \right\rceil$, 
the spherical zones are numbered from 1 to $q$. 

Now, we consider the following two steps:
\begin{itemize}
	\item a suitable family of $q$ spherical zones of equal width (that is, the width of the strip of the spherical zone) $\delta_{s_1} \equiv \delta_Z$ (possibly except for one of them) on the sphere and parallel to the $xy$-plane is constructed;
	\item the set ${\cal S}_n$ of nodes is partitioned by applying the spherical zone structure into $q$ subsets ${\cal S}_{kn_k}$, whose elements are $(x_{k1},y_{k1},z_{k1})$, $(x_{k2},y_{k2},z_{k2})$, $\ldots$, $(x_{kn_k},y_{kn_k},z_{kn_k})$, $k=1,2,\ldots,q$. 
\end{itemize}

\mascotsubsection{Localization stage}
In the searching procedure we consider three zones fo each node to be  examined.(for details see \cite{Cavoretto-DeRossi10}). Then, for each zone $k$, $k=1,2,\ldots,q$, a spherical zone searching routine is considered, examining the nodes from zone $k-i^*$ to zone $k+i^*$. (Note that if $k-i^* < 1$ or $k+i^* > q$ it will assign $k-i^* = 1$ and $k+i^*=q$, respectively.)

After defining the spherical zones to be examined for each node of ${\cal S}_{kn_k}$, $k=1,2,\ldots,q$, a spherical zone searching procedure is applied to determine all nodes belonging to a local neighbourhood. Here, we check whether the number of nodes in each neighbourhood is greater or equal to $n_Z$; if the condition is not satisfied, we repeat the process increasing the value of $k_1$ in (\ref{delta_Z}).

In the following, all the nodes belonging to a circular neighbourhood centred at ($x_i,y_i,z_i$), $i=1,2,\ldots,n$, are first ordered by applying a quicksort procedure on the distance, and then reduced to $n_Z$. Thus, taking only the $n_Z$ nodes closest to the centre ($x_i,y_i,z_i$), $i=1,2,\ldots,n$, of the neighbourhood, a local interpolant $Z_j$, $j=1,2, \ldots, n$, of the form given by (\ref{zonbasfun}), is constructed for each node.

\begin{figure}[ht!]
\begin{center}
  \includegraphics[width=8.cm]{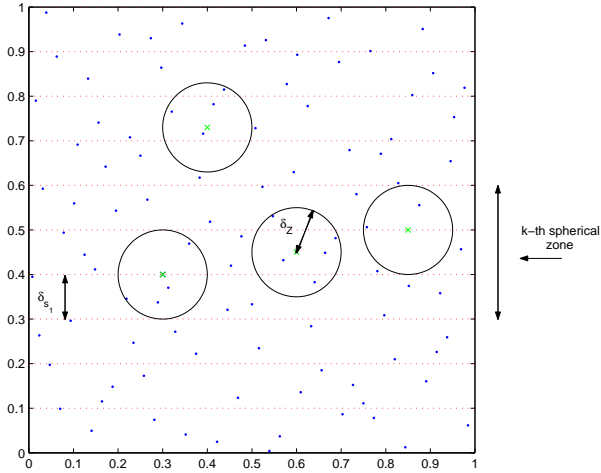}
\end{center}  
  \caption{Example of circular neighbourhoods.}
\label{mags}
\end{figure}

\mascotsubsection{Evaluation stage}
For each evaluation point $(x,y,z) \in {\cal E}_s$, a circular neighbourhood is constructed, whose spherical radius depends on the node number $n$, the parameter value $n_W$, and the (positive integer) number $k_2$, that is,
\begin{eqnarray}
\label{sphdelta_W}
	\delta_W = \arccos \left(1-2\sqrt{k_2} \frac{n_W}{n}\right), \hspace{0.5cm} k_2 = 1, 2, \ldots
\end{eqnarray}
After the number of spherical zones is determined by 
$r = \left\lceil \pi /\delta_W \right\rceil$, 
the spherical zones are numbered from 1 to $r$.

Then, a second family of $r$ spherical zones of equal width $\delta_{s_2} \equiv \delta_W$ (possibly except for one of them) on the sphere and parallel to the $xy$-plane is constructed. The sets ${\cal S}_n$ and ${\cal E}_s$ are partitioned into $r$ subsets ${\cal S}_{p_k}$ and ${\cal E}_{q_k}$, $k=1,2, \ldots, r$, respectively, so that the nodes of ${\cal S}_{p_k}$ and the evaluation points of ${\cal E}_{q_k}$ belong to the $k$-th zone.

A spherical zone searching procedure on the sphere is applied for each evaluation points of ${\cal E}_{p_k}$, $k=1,2,\ldots,r$, in order to find all nodes belonging to a (local) neighbourhood of centre $(x_i,y_i,z_i)$ and geodesic radius $\delta_W$. In this phase the basic idea is similar to that presented above, substituting $\delta_{s_1}$ by $\delta_{s_2}$, and $\delta_Z$ by $\delta_W$.

The nodes of each neighbourhood are first ordered by applying a quicksort procedure, and then reduced to $n_W$. Thus, considering only the $n_W$ nodes closest to the evaluation point $(x,y,z)$, a local weight function $\bar{W}_j(x,y,z)$, $j=1,2,\ldots,n$, is found. Finally, applying the formula (\ref{modsphshep}), the surface can be approximated at evaluation points $(x,y,z) \in {\cal E}_s$.


\mascotsection{Numerical experiments}
In this section we test accuracy of the local Shepard's method, which makes use of ZBFs plus spherical harmonics. In doing so, we take three scattered data sets ($n=1000, 4000, 16000$) obtained by using a MATLAB code (see \cite{Fornberg-Piret07,Cavoretto-DeRossi12}), which generates uniformly random node distributions on the sphere. The interpolant is evaluated on a set of $s=600$ spiral points generated by the algoritm of Saff and Kuijlaars (see \cite{Cavoretto-DeRossi10}), which gives us a fairly good point distribution over $\mathbb{S}^{2}$, tracing out an imaginary spiral from the South pole to the North pole.


Data values are taken from the restriction of the following trivariate test functions on $\mathbb{S}^{2}$, that is,
\begin{equation}
f_1(x,y,z)= {{\rm e}^x+2{\rm e}^{y+z}}/{ 10}, \hspace{0.3cm} f_2(x,y,z)= \sin x\ \sin y\ \sin z.
\end{equation}

Before analyzing numerical experiments obtained using the local interpolation scheme, we remark that the choice of the localization parameters $n_Z$ and $n_W$ is of great importance, because it determines the level of accuracy of the method and the  efficiency of the corresponding algorithm. Anyway, in our tests we consider a good trade-off beetween these two remarkable items, taking $n_Z=15$ and $n_W=10$. Similarly, we act to choose the value of ZBF shape parameter; indeed, we pick the mean value of the interval $(0,1)$, i.e. assuming $\gamma = 0.5$. Specifically, this fact leads to a good comprimise between accuracy and stability.   

Moreover, since the space $H_d$ of the spherical harmonics of degree $d$ has dimension $N_{d,3}= 2 d + 1$ on $\mathbb{S}^{2}$, and fixing $L \geq 0$ as the desired degree of the spherical harmonic component of the approximation, we have that $U={\rm dim} {\cal H}_L= (L+1)^2$. It follows that the necessary condition $n_Z \geq (L+1)^2$ imposes an upper limit of $L=2$. Thus, in the numerical experiments we take $L=-1,0,1,2$, where $L=-1$ denotes no spherical harmonic component.

Then, we compute the relative root mean square errors (RRMSEs) on the evaluation points, 
reporting in Tables \ref{tab1_errors} -- \ref{tab3_errors} the achieved results.

\begin{figure}[ht!]
\begin{center}
  \includegraphics[width=6.5cm]{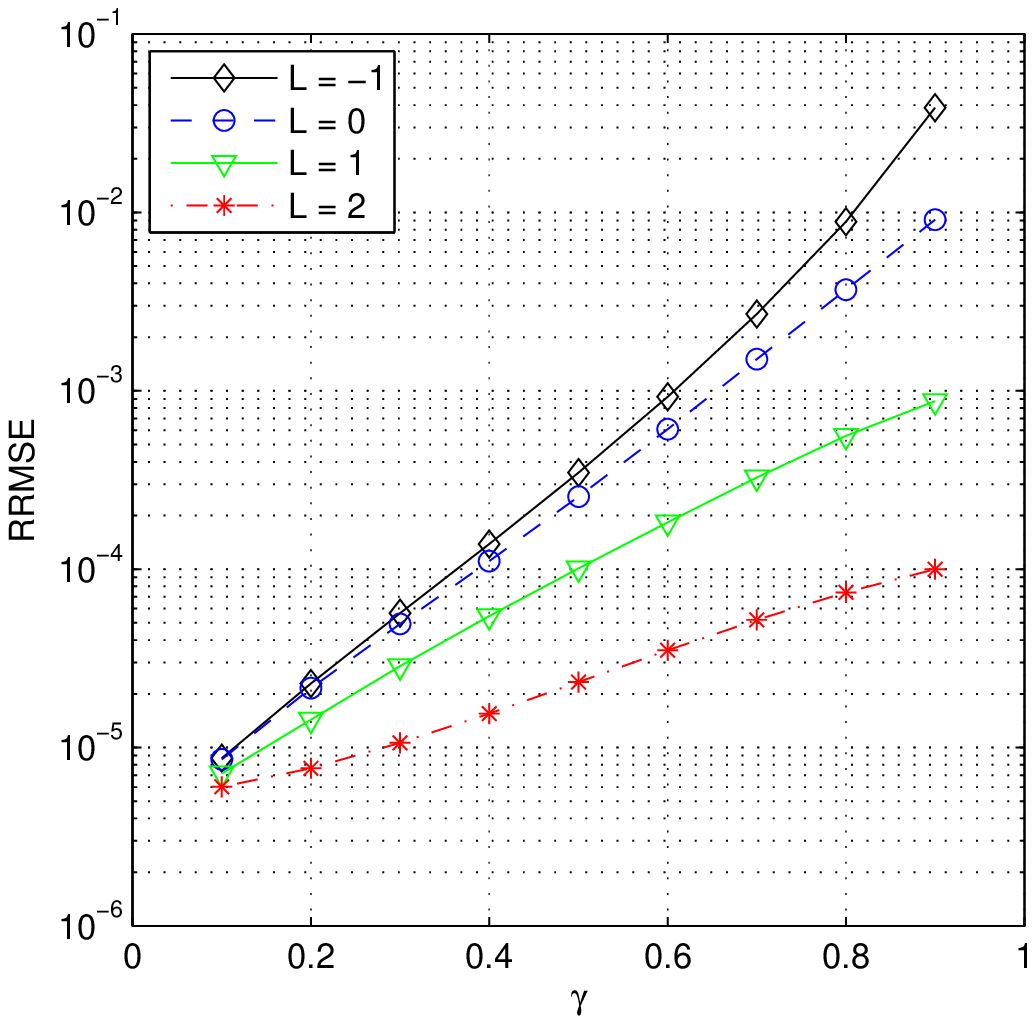} \hskip -1cm
  \includegraphics[width=6.5cm]{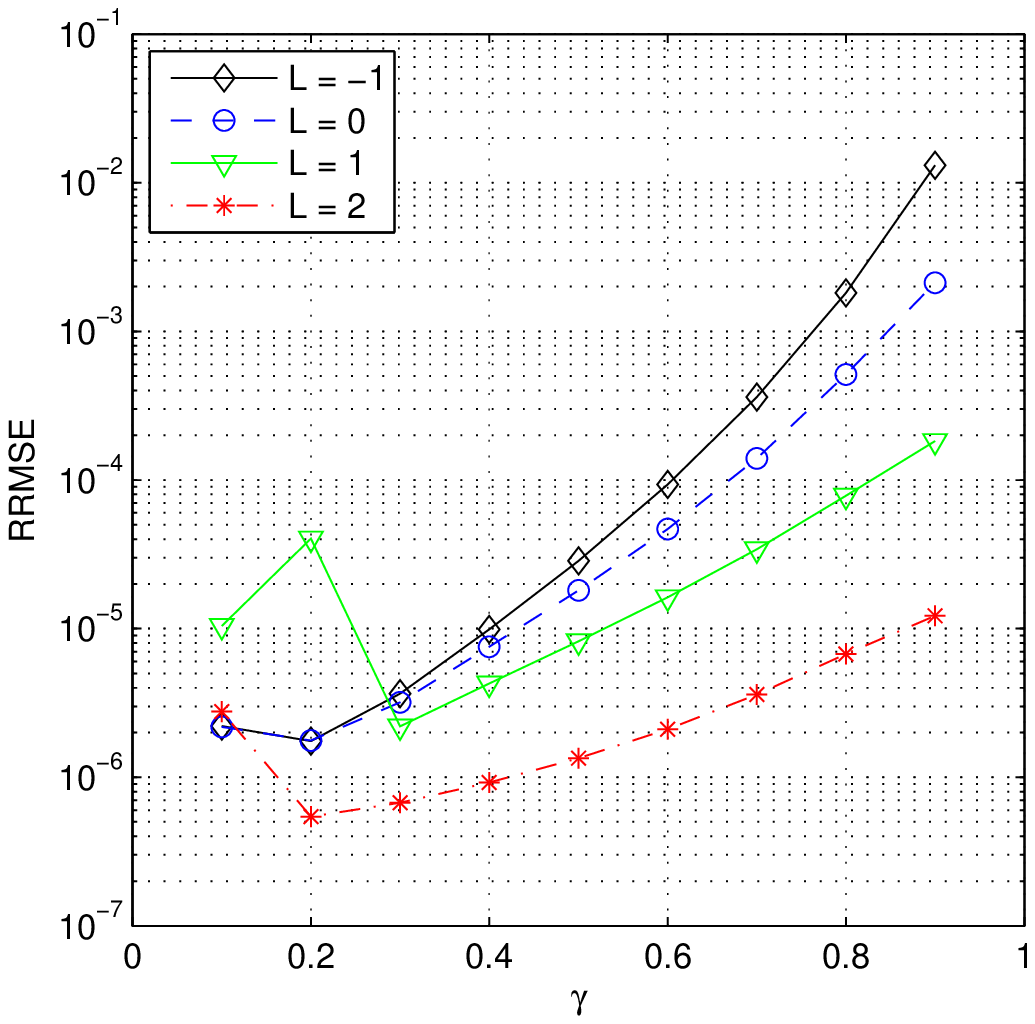} \hskip -1cm
\end{center}  
  \caption{RRMSEs with $n =1000$ (left) and $n = 4000$ (right) by varying $\gamma$ for $f_1$.}
\label{M_shape_par}
\end{figure}

Even though interpolation errors are rather low also when there is no spherical harmonic component (i.e., $L = -1$), these results point out that the addition of a spherical harmonics may produce an improvement of accuracy, mainly when the value of $L$ increases. This effect is noted for each of the considered node sets. In Figure \ref{M_shape_par} we show the behaviour of the RRMSEs by varying the IMQ shape parameter $\gamma$ in the interval $(0,1)$ for $f_1$.

\begin{table}[ht!]
		\begin{center}
			\begin{tabular}{cccc} 
			 \hline
			 \rule[0mm]{0mm}{3ex}
			  	$L \setminus n$     & $1000$ & $4000$ & $16000$ \\
				\hline
				\rule[0mm]{0mm}{3ex}
				  $-1$ \ & \ $3.4759{\rm e}-4$  \ & \ $2.8568{\rm e}-5$ \ & \ $1.7244{\rm e}-6$      \\
			  \rule[0mm]{0mm}{3ex}
			    $0$  \ & \ $2.5466{\rm e}-4$  \ & \ $1.8057{\rm e}-5$ \ & \ $1.2770{\rm e}-6$      \\
			  \rule[0mm]{0mm}{3ex}
				  $1$  \ & \ $1.0109{\rm e}-4$  \ & \ $8.2052{\rm e}-6$ \ & \ $8.1097{\rm e}-7$      \\
				\rule[0mm]{0mm}{3ex}
				  $2$  \ & \ $2.3277{\rm e}-5$  \ & \ $1.3413{\rm e}-6$ \ & \ $4.3374{\rm e}-8$      \\
			  \hline 
			\end{tabular}
		\end{center}
			\caption{RRMSEs obtained by IMQ with $\gamma = 0.5$ for $f_1$.}
			\label{tab1_errors}
	\end{table}

\begin{table}[ht!]
		\begin{center}
			\begin{tabular}{cccc} 
			 \hline
			 \rule[0mm]{0mm}{3ex}
			  	$L \setminus n$     & $1000$ & $4000$ & $16000$ \\
				\hline
				\rule[0mm]{0mm}{3ex}
				  $-1$ \ & \ $2.6059{\rm e}-2$  \ & \ $5.5551{\rm e}-3$ \ & \ $4.2012{\rm e}-5$      \\
			  \rule[0mm]{0mm}{3ex}
			    $0$  \ & \ $2.5769{\rm e}-2$  \ & \ $5.6371{\rm e}-3$ \ & \ $4.2514{\rm e}-5$      \\
			  \rule[0mm]{0mm}{3ex}
				  $1$  \ & \ $3.9581{\rm e}-2$  \ & \ $6.1304{\rm e}-3$ \ & \ $6.1078{\rm e}-5$      \\
				\rule[0mm]{0mm}{3ex}
				  $2$  \ & \ $6.9575{\rm e}-3$  \ & \ $3.4626{\rm e}-4$ \ & \ $1.0221{\rm e}-5$      \\			  
			  \hline 
			\end{tabular}
		\end{center}
			\caption{RRMSEs obtained by IMQ with $\gamma = 0.5$ for $f_2$.}
			\label{tab3_errors}
	\end{table}

Finally, we also apply our local interpolation scheme to geomagnetic data, known as \textsl{MAGSAT} (MAGnetic field SATellite) \cite{MAGSAT}. In particular, here we consider two subsets of nodes ($n=2084, 4088$) obtained after manipulating and refining the original MAGSAT data, so that the distribution of each set is reasonably uniform on $\mathbb{S}^{2}$. Specifically, we randomly select from the original data sets $n$ geomagnetic nodes for the interpolation process, taking $s=200$ points for the cross-validation.  As an example, the representation of the $n=2084$ nodes  is shown in Figure \ref{mags}. Then, in Table \ref{tab1} we report RRMSEs obtained by using MAGSAT data, taking $\gamma=0.96$, $n_Z=12$ and $n_W=10$ as parameters, for $L=-1,0$.

\begin{figure}[ht!]
\begin{center}
  \includegraphics[width=9.cm]{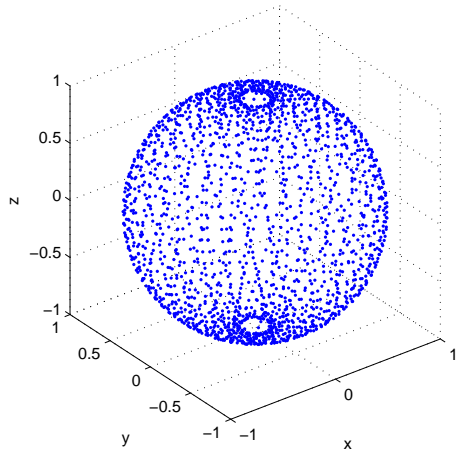}
\end{center}  
  \caption{Example of MAGSAT data.}
\label{mags}
\end{figure}

\begin{table}[htbp]
		\begin{center}
			\begin{tabular}{|c|c|c|c|c|} \hline
$n$& \multicolumn{2}{c|}{  \rule[-2mm]{0mm}{6mm}  2084} & \multicolumn{2}{c|}{  \rule[-2mm]{0mm}{6mm} 4088} \\
			  \cline{2-5} \rule[-2mm]{0mm}{6mm}
			  $L$	& $-1$ & $0$ & $-1$ & $0$ \\
				\hline
				\rule[0mm]{0mm}{3ex}
				IMQ & $4.6865{\rm e}-2$ & $2.2349{\rm e}-2$   & $4.1185{\rm e}-2$ & $2.3109{\rm e}-2$ \\
			  \hline 
			\end{tabular}
		\end{center}
			\caption{RRMSEs for MAGSAT data.}
		\label{tab1}
	\end{table}


\subsection{Acknowledgements}
The author gratefully acknowledge the financial support of the GNCS-INDAM.


\end{document}